\documentclass[11pt]{article}
\def\boe{\begin{enumerate}}
\def\eoe{\end{enumerate}}

\newtheorem{proposition}{{\bf Proposition}}
\newtheorem{lemma}{{\bf Lemma}}
\newtheorem{corollary}{{\bf Corollary}}
\newtheorem{assumption}{{\bf Assumption}}
\newtheorem{definition}{{\bf Definition}}
\newtheorem{theorem}{{\bf Theorem}}

\newcommand\ca[1]{{\cal{#1}}}
\newcommand\lo[1]{_{\nano{#1}}}
\newcommand\hi[1]{^{\nano{#1}}}

\def\proof{\noindent {\sc Proof. }}

\def\L{{\cal L}}

\def\trans{^{\mbox{\tiny{\sf  T}}}}

\def\inv{^{\mbox{\tiny $-1$}}}

\def\var{\mathrm{var}}
\def\cov{\mathrm{cov}}

\newcommand{\indep}{\;\, \rule[0em]{.03em}{.65em} \hspace{-.41em}
\rule[-.02em]{.65em}{.03em} \hspace{-.41em}
\rule[0em]{.03em}{.65em}\;\,}

\def\eop{\hfill $\Box$ \\
}
\def\trans{^{\mbox{\tiny{\sf T}}}}

\def\iff{\Leftrightarrow}

\def\ali{&\,}

\def\ran{\mathrm{ran}}

\def\of{{\nano {\circ}}}
\def\nano{\scriptscriptstyle}

\def\real{{\mathbb R}}

\def\ka{\kappa}
\def\cenclass{\mathfrak{S} \lo {Y|X}}

\def\oc{^{\nano \perp}}

\def\L2T{L \lo 2 (T)}
\def\L2TX{L \lo 2 (T\lo X)}
\def\L2TX{L \lo 2 (T\lo Y)}

\def\ali{&\,}
\def\spn{{\rm{span}}}
\def\ali{&\,}

\def\eod{\end{document}}

\def\cran{{\overline{{\rm{ran}}}}}

\newfont{\rsfsten}{rsfs10 scaled 1100}
\newfont{\rsfstena}{rsfs10 scaled 800}
\newfont{\rsfstenb}{rsfs10 scaled 500}

\def\iff{\Leftrightarrow}

\usepackage{multirow,verbatim}
\usepackage{latexsym, amsbsy, amssymb}
\usepackage{natbib}
\usepackage{tikz}
\usepackage{graphicx}
\usepackage{epsfig}
\usepackage{amsmath}
\usepackage{multirow}
\usepackage{verbatim}
\usepackage{color}

\def\var{{\mathrm{var}}}

\def\eop{\hfill $\Box$ \\ }
\def\proof{\noindent{\sc{Proof. \ }}}

\def\nano{\scriptscriptstyle}
\def\real{\mathbb R}

\def\var{\mathrm{var}}

\def\oc{\hi{\perp}}

\def\ran{\mathrm{ran}}
\def\ker{\mathrm{ker}}
\def\cran{\overline{\mathrm{ran}}}
\def\iff{\Leftrightarrow}
\def\cov{\mathrm{cov}}
\def\nano{\scriptscriptstyle}
\def\inv{\hi{\nano -1}}

\def\nano{\scriptscriptstyle}
\def\of{\mbox{\raisebox{1pt}{$\nano{\circ}$}}}
\def\ka{\kappa}

\def\ali{&\,}

\def\hii#1{\hi{(#1)}}

\def\spn{\mathrm{span}}

\newcommand{\RN}[1]{%
  \textup{\uppercase\expandafter{\romannumeral#1}}%
}

\def\ca#1{{\cal#1}}
\def\lo{_}
\def\hi{^}
\def\ran{\mathrm{ran}}
\def\cran{\overline{\ran}}
\def\spn{\mathrm{span}}
\def\cspn{\overline{\spn}}
\def\cenclass{\mathfrak{S} \lo {Y|X}}

\begin{document}

\title{On relative universality, regression operator,  and conditional independence }

\author{Bing Li, Ben Jones and Andreas Artemiou}

\maketitle

\begin{abstract}
 The notion of relative universality with respect to a $\sigma$-field was introduced to establish the unbiasedness and Fisher consistency of an estimator in nonlinear sufficient dimension reduction. However, there is a gap in the proof of this result in the existing literature. The existing definition of relative universality seems to be   too strong for the proof to be valid. In this note we modify the definition of relative universality using the concept of $\epsilon$-measurability, and rigorously establish the mentioned unbiasedness and Fisher consistency. The significance of this result is beyond its original context of sufficient dimension reduction, because relative universality allows us to use the regression operator to fully characterize conditional independence, a crucially important statistical relation that sits at the core of  many areas and methodologies in  statistics and machine learning, such as dimension reduction, graphical models, probability embedding, causal inference, and Bayesian estimation.
\end{abstract}


\section{Introduction}

In this paper we rigorously introduce the notion of relative universality, a  critical assumption needed for characterizing conditional independence. We then use this concept to rigorously establish a relation between the regression operator \citep{fukumizu2007statistical,lee2016} and conditional independence. Relative universality was first introduced in \citet[][Section 13.4]{li2018sufficient} as a mechanism to establish the unbiasedness and Fisher consistency of  Generalized Sliced Inverse Regression (GSIR), an important estimator  for  nonlinear sufficient dimension reduction. Nonlinear sufficient dimension reduction is a methodology to reduce the dimension of the high-dimension predictor in a regression setting, and has undergone vigorous development during the recent years. See, for example,   \citet{wu2008}, \citet{yeh-huang-lee-2009}, \citet{li-artemiou-chiaromonte-2011}, \citet{lee2013},   \cite{li2017}, and \cite{zhang2024nonlinear}.
However, two coauthors (second and third) of the current paper have recently discovered   a gap in the proof of the above result if we use the original definition of relative universality in \cite{li2018sufficient}. The goal of this paper is, first, to correct the error in \cite{li2018sufficient}, and second,  to systematically  and rigorously develop the theory surrounding relative universality, regression operator, and conditional independence. Since conditional independence is a widely used mechanism in many areas of  statistics and machine learning, such as sufficient dimension reduction \citep{Li1991, Cook1994, li2018sufficient}, sufficient graphical models \citep{li2024sufficient}, nonparametric variable selections \citep{lee2016},   causal estimation, and Bayesian inference \citep{li2019graduate},  a carefully and systematically developed theory of relative universality would be conducive for developing methodologies surrounding conditional independence in various theoretical and applied settings.

The rest of the paper is organized as follows. In Section \ref{section:relative universality}, we review the notion of relative universality as originally   defined by \cite{li2018sufficient}, explain why it is unfit to prove the unbiasedness and Fisher consistency of GSIR, and then relax it using $\epsilon$-measurability. We prove a theorem that makes the modified relative universality useful.  In Section \ref{section:regression operator}, we introduce nonlinear sufficient dimension reduction, the  regression operator, and the generalized Sliced Inverse Regression.   In Section \ref{section:unbiased}, we establish the main result of this paper --- how the regression operator characterizes conditional independence --- using the modified version of relative universality. We will also describe the gap in the proof in \cite{li2018sufficient} that motivates this paper.  In Section \ref{section:conclusion}, we summarize the main message of this paper and give an overview of the logic line underlying our  development.



\section{Relative universality}\label{section:relative universality}

The concept  of universality was introduced by \cite{micchlli2006universal} to describe the richness of the reproducing kernel Hilbert space (RKHS) generated by a positive definite kernel: we say that a kernel is universal if the RKHS it generates is dense in the class of bounded and continuous functions. See also \cite{Sriperumbudur2011universality}. It is widely used in the development of RKHS-related methodologies. See, for example, \cite{Caponnetto2008kernel}, \cite{fukumizu2009}, and \cite{simon2018kernel}.
As a device to handle conditional independence, \cite{li2018sufficient} introduced the notion of relative universality, which is, loosely,  universality with respect to a $\sigma$-field --- the $\sigma$-field being conditioned on in the conditional independence. In the following we first describe this concept, why it is unfit to prove the unbiased of GSIR,  and then relax it so as to facilitate the proof.

\subsection{\cite{li2018sufficient}'s definition of relative universality}

To motivate our development, we first review the definition of relative universality given in \cite{li2018sufficient}, Section 13.4. Since our development will be more general than the RKHS framework of \cite{li2018sufficient}, we will state that definition in the more general setting.

Let $(\Omega, \ca F, P)$ be a probability space. Let $(\Omega \lo X, \ca F \lo X)$ a measurable space. Let $X: \Omega \to \Omega \lo X$ be a random vector, and $P \lo X = P \of X \inv$ its distribution. Let $L \lo 2 (P \lo X)$ be the set of all square-integrable functions on $\Omega \lo X$ with respect to $P \lo X$. Without loss of generality, we assume that $\ca F$ is the $\sigma$-field generated by $X$; that is, $\ca F = X \inv (\ca F \lo X)$.  Let $\ca H \lo X \subseteq L \lo 2 (P \lo X)$ be a Hilbert space. Note that the inner product in $\ca H \lo X $ need not be the same as that in $L \lo 2 (P \lo X)$, but we will make the following assumption.
\begin{assumption} \label{assumption:norm norm}  There is a constant $C > 0$ such that
\begin{align}\label{eq:norm norm}
\| f \| \lo {L \lo 2 (P \lo X)} \le C \| f \| \lo {\ca H \lo X }.
\end{align}
\end{assumption}

The scenario that Li (2018) considered is where   $\ca H \lo X $ is the RKHS generated by $\ka \lo X : \Omega \lo X \times \Omega \lo X \to \real$ where $E [\ka \lo X (X,X)] < \infty$. This Hilbert space does satisfy Assumption \ref{assumption:norm norm} because, for any $f \in \ca H \lo X $,
\begin{align*}
\| f \| \lo {L \lo 2 (P \lo X)} \hi 2 = E \left( \langle \ka \lo X (\cdot, X), f \rangle \hi 2 \lo {\ca H \lo X} \right)  \le \| f \| \lo {\ca H \lo X} \hi 2 \, E \ka \lo X (X, X).
\end{align*}
So  (\ref{eq:norm norm}) is satisfied with $C = \sqrt{E[\ka \lo X (X,X)]}$.

\def\iff{\Leftrightarrow}

In the following, for a subset $\ca S$ of $L \lo 2 (P \lo X)$ and a sub-$\sigma$-field $\ca G$ of $\ca F$, let
$\ca S \lo {\ca G}$ be the collection of functions in $\ca S$  {such that, for any $f \in \ca S$, $f ( X)$ is } measurable with respect to $\ca G$. Thus, for example, $(\ca H \lo X) \lo {\ca G}$ is the set of all functions $f$ in $\ca H \lo X $  {such that $f (X)$ is } measurable with respect to $\ca G$, and $L \lo 2 (P \lo X) \lo {\ca G}$ is the set of all functions  {$f$} in $L \lo 2 (P \lo X)$  {such that $f(  X)$ is } measurable with respect to $\ca G$. The next proposition shows that $\ca S \lo {\ca G}$ is a closed linear subspace of $\ca S$ if $\ca S$ is a Hilbert space satisfying (\ref{eq:norm norm}). This result was assumed in \cite{li2018sufficient}  without proof.

\begin{proposition}  {If $\ca H \lo X \subseteq L \lo 2 (P \lo X)$ is a Hilbert space satisfying Assumption \ref{assumption:norm norm},} and $\ca G$ is a sub $\sigma$-field of $\ca F$, then $(\ca H \lo X) \lo {\ca G}$ is a closed linear subspace of $\ca H \lo X $.
\end{proposition}

\proof That $(\ca H \lo X ) \lo {\ca G}$ is a linear subspace is obvious. We now prove it is closed.  {Let $h \in \ca H \lo X $ be an accumulation point of $(\ca H \lo X) \lo {\ca G}$. We need to show tht $h$ is a member of $(\ca H \lo X) \lo {\ca G}$. Let }  $\epsilon > 0$ and  $g$   a member of $(\ca H \lo X ) \lo {\ca G}$ satisfying $\| g - h \| \lo {\ca H \lo X  } < \epsilon$. Then,  {by Assumption \ref{assumption:norm norm}},
\begin{align*}
\var [ g(X) - h(X)] \le E \{  [ g(X) - h(X)] \hi 2\} \le  {C}\| g - h \| \lo {\ca H \lo X  } \hi 2 < C \epsilon \hi 2.
\end{align*}
It follows that
\begin{align*}
E \{ \var[ g(X) - h(X) | \ca G ] \} \le \var [ g(X) - h(X)]  < C \epsilon \hi 2.
\end{align*}
Since  {$g (X)$} is measurable with respect to $\ca G$, the left-hand side is  {$E \{ \var[ h(X) | \ca G ] \}$, and hence $E \{ \var[ h(X) | \ca G ] \}  < \epsilon \hi 2$}. Since $\epsilon$ can be an arbitrarily small constant, we have $E \{ \var[ h(X) | \ca G ] \}  =0$, implying  $\var[ h(X) | \ca G ] = 0$ almost surely. But this means $h(X)$ is a constant given $\ca G$. Thus $h(X)$ is measurable with respect to $\ca G$. \eop

We now give the formal definition of relative universality in \cite{li2018sufficient}. Since we are going to relax this condition, we will call it {\em strong} relative universality and save the term ``relative universality'' for the modified version.
\begin{definition}\label{definition:2018} For a given sub-$\sigma$-field $\ca G$ of $\ca F$, we say that $\ca H \lo X$ is strongly relatively universal with respect to $\ca G$ if $(\ca H \lo X ) \lo {\ca G}$ is dense in $L \lo 2 (P \lo X) \lo {\ca G}$ modulo constants.
\end{definition}


\subsection{Relaxation of \cite{li2018sufficient}'s definition of relative universality}

To relax Definition \ref{definition:2018}, we first introduce $\epsilon$-measurability.

\begin{definition} For a given $\epsilon \ge  0$, we
      say that $f(X)$ is $\epsilon$-measurable with respect to $\ca G$ if
      \begin{align}
          E(\var[f(X)|\ca G]) < \epsilon.
      \end{align}
\end{definition}
This is a generalization of measurability: if $\epsilon = 0$, then $\epsilon$-measurable means $E(\var[f(X)|\ca G]) = 0$, which implies $\var[f(X)|\ca G] = 0$ almost surely, which implies $f(X)$ is measurable $\ca G$.
For any $\epsilon > 0$,  let
$
(\ca H \lo X ) \lo {\ca G} (\epsilon)$ be the collection of all $\epsilon$-measurable functions in $\ca H \lo X$. We have the following proposition.

\begin{proposition} If $(\ca H \lo X ) \lo {\ca G}$ and $(\ca H \lo X ) \lo {\ca G} (\epsilon)$ are as defined above, then
\begin{align*}
(\ca H \lo X ) \lo {\ca G} = \cap \lo {\epsilon > 0} (\ca H \lo X ) \lo {\ca G} (\epsilon).
\end{align*}
\end{proposition}

\proof Note that
\begin{align*}
(\ca H \lo X ) \lo {\ca G} =\ali  \{ f \in \ca H \lo X : \mbox{$f(X)$ is measurable with respect to $\ca G$} \} \\
= \ali \{ f \in \ca H \lo X : E(\var[f(X)|\ca G]) = 0 \} \\
= \ali \{ f \in \ca H \lo X : E(\var[f(X)|\ca G]) < \epsilon  \ \mbox{for all $\epsilon > 0$} \}.
\end{align*}
where the right-hand side is, by definition, $\cap \lo {\epsilon > 0} (\ca H \lo X ) \lo {\ca G} (\epsilon)$.
\eop

Note that $(\ca H \lo X ) \lo {\ca G} (\epsilon)$ increases with $\epsilon$ in the sense that, if $0 \le \epsilon \lo 1< \epsilon \lo 2$, then $(\ca H \lo X ) \lo {\ca G} (\epsilon \lo 1) \subseteq (\ca H \lo X ) \lo {\ca G} (\epsilon \lo 2)$.
We now introduce our new definition of relative universality.


\begin{definition} We say that $\ca H \lo X $ is relatively universal with respect to $\ca G$ if, for any $\epsilon > 0$,
$ (\ca H \lo X ) \lo {\ca G} (\epsilon)$ is dense in
$ L \lo 2 (P \lo X) \lo {\ca G} $ modulo constants.
\end{definition}

The difference between the new definition and the original definition in Li (2018) is that we replaced measurable functions by $\epsilon$-measurable functions for any $\epsilon > 0$. Thus  this condition is weaker than Li (2018) definition. Nevertheless, the two definitions are very close because, when $\epsilon$ is small, a function in $(\ca H \lo X ) \lo {\ca G}(\epsilon)$ is very nearly measurable with respect to $\ca G$.

For further development, it is useful to restate the above  definition in an alternative, but equivalent form. There will be three types of orthogonality involved in our discussion. We denote the orthogonality in $\ca H \lo X $ by $\perp \lo 1$, the orthogonality in $L \lo 2 (P \lo X)$ by $\perp \lo 2$, and the orthogonality in $L \lo 2 (P \lo X)$ modulo constant by $\perp \lo 3$. That is:
\begin{enumerate}
\item For $f, g \in \ca H \lo X $,
$f \perp \lo 1 g \iff \langle f, g \rangle \lo {\ca H \lo X }=0$;
\vspace{-0.08in}
\item For $f, g \in L \lo 2 (P \lo X)$, $f \perp \lo 2 g \iff  {E[f(X), g(X)] = 0}$;
\vspace{-0.08in}
\item For $f, g \in L \lo 2 (P \lo X)$,  $f \perp \lo 3 g \iff \cov[f(X), g(X)] = 0$.
\end{enumerate}
Orthogonal complements are defined accordingly: for example, for a set $A \subseteq L \lo 2 (P \lo X)$, $A \hi {\perp \lo 3}$ the set
\begin{align*}
\{ f \in L \lo 2 (P \lo X): \cov[f(X), g(X)] = 0 \quad \mbox{for all $g \in A$}\}.
\end{align*}
It is well known that, for a generic Hilbert space $\ca H \lo X $ and its subsets $A$ and $B$ with $A \subseteq B$, $A$ is dense in $B$ if and only if $A \hi \perp = B \hi \perp$. This statement also holds for $\perp \lo 3$ if we replace ``dense'' with ``dense modulo constants''.

\begin{lemma}\label{lemma:perp 3} Suppose $A$ and $B$ are subsets of $L \lo 2 (P \lo X)$ with $A \subseteq B$. Then $A$ is dense in $B$ modulo constants if and only if $A \hi {\perp \lo 3} = B \hi {\perp \lo 3}$.
\end{lemma}

\proof Introduce the following equivalence relation in $L \lo 2 (P \lo X)$:
\begin{align*}
f \sim g \quad \iff \quad f(X) - g(X) = \mbox{constant almost surely}.
\end{align*}
Then the quotient space $L \lo 2 (P \lo X) /\!\! \sim$, equipped with the inner product  $ \cov[f(X), g(X)]$, forms a Hilbert space. The result then follows from Corollary 1.10 of Conway (1990). \eop

Using this result we immediately arrive at the following equivalent conditions of relative universality.

\begin{corollary} The following statements are equivalent:
\begin{enumerate}
\item for every $\epsilon > 0$, $(\ca H \lo X ) \lo {\ca G} (\epsilon) $ is dense in $ {L \lo 2(P \lo X) \lo {\ca G}}$ modulo constant;
  \vspace{-0.08in}
\item  for every $\epsilon > 0$, $[(\ca H \lo X ) \lo {\ca G}  (\epsilon)] \hi {\perp \lo 3} \subseteq [L \lo 2 (P \lo X) \lo {\ca G} ]\hi {\perp \lo 3}$; \vspace{-0.08in}
\item  for every $\epsilon > 0$, $[(\ca H \lo X ) \lo {\ca G} (\epsilon) ] \hi {\perp \lo 3} = [L \lo 2 (P \lo X) \lo {\ca G} ]\hi {\perp \lo 3}$.
\end{enumerate}
\end{corollary}

\proof The equivalence of {\em 1.} and {\em 3.} follows from Lemma \ref{lemma:perp 3}; the equivalence of {\em 2.} and {\em 3.} follows from
$[(\ca H \lo X ) \lo {\ca G} ] \hi {\perp \lo 3} \supseteq [L \lo 2 (P \lo X) \lo {\ca G} ]\hi {\perp \lo 3}$, which is obviously true.  \eop


To gain more intuition about this concept, it is helpful to consider the special cases where $\ca G$ is the largest $\sigma$-field $\sigma (X)$ and the smallest $\sigma$-field $\{\varnothing, \Omega \}$.

\begin{corollary}\label{corollary:two special cases}  {If $\ca H \lo X \subseteq L \lo 2 (P \lo X)$ is a Hilbert space}, then the following statements hold true.
\begin{enumerate}
\item $\ca H \lo X $ is relatively universal with respect to $\sigma (X)$ if and only if $\ca H \lo X $ is dense in $L \lo 2 (P \lo X)$ modulo constant; \vspace{-0.08in}
\item $\ca H \lo X $ is always relatively universal with respect to $\{\varnothing, \Omega \}$.
\end{enumerate}
\end{corollary}


\def\oc3{\hi {\perp \lo 3}}

\proof 1. Note that, for any $\epsilon > 0$,
\begin{align}\label{eq:largest sigma field result}
(\ca H \lo X )\lo {\sigma (X)} (\epsilon) = \{ f \in \ca H \lo X  : E(  \var[f(X)| X] ) < \epsilon \}= \{ f \in \ca H \lo X  : 0 < \epsilon \}=\ca H \lo X.
\end{align}
 {Also note that $L \lo 2 (P \lo X) \lo {\sigma (X)} = L \lo 2 (P \lo X)$}. If $(\ca H \lo X ) \lo {\ca G} (\epsilon)  \hi {\perp \lo 3} \subseteq [L \lo 2 (P \lo X) \lo {\ca G} ]\hi {\perp \lo 3}$ for any $\epsilon > 0$, then $\ca H \lo X   \hi {\perp \lo 3} \subseteq [L \lo 2 (P \lo X) \lo {\ca G} ]\hi {\perp \lo 3}$, which is equivalent to saying $\ca H \lo X $ is dense in $L \lo 2 (P \lo X)$ modulo constant.  {Conversely, if $\ca H \lo X$ is dense in $L \lo 2 (P \lo X)$ modulo constant, then $\ca H \lo X \oc3 = L \lo 2 (P \lo X) \oc3$ which, by (\ref{eq:largest sigma field result}), implies $[( \ca H \lo X ) \lo {\sigma (X)} (\epsilon) ] \oc3 \subseteq L \lo 2 (P \lo X) \lo {\sigma (Y)}$. Thus $\ca H \lo X$ is relatively universal with respect to $\sigma (X)$. }

\noindent 2.
 {If $\ca G = \{\varnothing, \Omega\}$, then $L \lo 2 (P \lo X) \lo {\ca G}$ is simply the set of all constant functions; that is, $f(x) = c$ almost surely for $c \in \real$. Hence $L \lo 2 (P \lo X) \lo {\ca G} \oc3 = L \lo 2 (  P \lo X)$, which contains the set $(\ca H \lo X) \lo {\ca G} \oc3$.} \eop



The next theorem is the fundamental property of relative universality that makes the concept useful. In fact, the   unbiasedness proof of Li (2018) is motivated by it.

\begin{theorem}\label{theorem:relative universality}  Given any sub-$\sigma$-field $\ca G$ of $\ca F$, if $ \ca H \lo X $ is dense in  $\ca G$ modulo constants,  then $\ca H \lo X $ is relatively universal with respect to $\ca G$.
\end{theorem}

\proof Since  $\ca H \lo X $ is relatively universal with respect to $\ca F$, by Corollary \ref{corollary:two special cases}, it is dense in $L \lo 2 (P \lo X)$ modulo constant. Let $\epsilon > 0$ and $\ca G$ be a sub $\sigma$-field of $\ca F$. Let $f $ be a member of $L \lo 2 (P \lo X)$ such that  $f \perp \lo 3 (\ca H \lo X ) \lo {\ca G} (\epsilon)$. Let $h \in L \lo 2 (P \lo X) \lo {\ca G}$. Let $\eta$ be a number such that $0 < \eta < \epsilon$, and  let $g \in \ca H \lo X $ be such that  $\var[h(X)-g(X)] < \eta$. Then
\begin{align*}
E \{ \var[g(X)| \ca G] \}=E \{ \var[h(X)-g(X)| \ca G] \}\le \var[h(X)-g(X)] < \eta
\end{align*}
Hence $g \in (\ca H \lo X ) \lo {\ca G} (\epsilon)$, and consequently $\cov[f(X), g(X)] = 0$. It follows that
\begin{align*}
\cov[f(X), h(X)]  \hi 2 = \ali \{  \cov[f(X), h(X)-g(X)] +  \cov[f(X), g(X)] \} \hi 2 \\
=\ali  \cov[f(X), h(X)-g(X)]\hi 2 \le \var[f(X)] \eta.
\end{align*}
Since the right-hand side can be arbitrarily small, we have $\cov[f(X), h(X)] =0$. Hence $f \perp \lo 3 L \lo 2 (P \lo X) \lo {\ca G}$. Thus we have shown
$(\ca H \lo X ) \lo {\ca G} (\epsilon) \hi {\perp \lo 3} \subseteq [L \lo 2 (P \lo X) \lo {\ca G}]  \hi {\perp \lo 3}$, as desired.  \eop


\section{Regression operator, nonlinear SDR and GSIR}\label{section:regression operator}

  {First, we} outline the construction of GSIR and some related terminologies. Our setting is more general than Lee, Li, and Chiaromonte (2013), Li (2018),  and  Li and Song (2017).

\subsection{Mathematical background and notations}

\def\cspn{\overline{\spn}}

For two   Hilbert spaces $\ca H$ and $\ca K$, the set of all bounded linear operators from $\ca H$ to $\ca K$  {is} written as $\ca B ( \ca H, \ca K)$. For a bounded linear operator $A \in \ca B ( \ca H, \ca K)$, we use $\ker (A)$ to denote the kernel of $A$: $\ker (A) = \{h \in \ca H: Ah = 0 \}$; we use $\ran (A)$ to denote the range of $A$: $\ran (A) = \{A h: h \in {\ca H} \}$. Recall that $\ker (A)$ is always a closed linear subspace, but $\ran (A)$ is a linear subspace that may not be closed. We use $\cran (A)$ to denote the closure of $\ran (A)$. For a subset  $\ca V$ of $\ca H$, we use $\spn ( \ca V)$ to denote the linear span of $\ca V$; that is,  the set of all finite linear combinations of members of $\ca V$. We use $\cspn (\ca V)$ to denote the closure of $\spn (\ca V)$.

In general, $A: \ca H \to \ca K$  {is a mapping} from $\ca H$ to $\ran (A)$, and this mapping may or may not be injective. But if we restrict $A$ on $\cran (A)$, then $(A | \,  \cran (A)): \cran(A) \to \ran(A)$ is always injective. Thus we can define a linear operator $A \hi \dagger : \ran (A) \to \cran(A)$ such that, for each $h \in \ran (A)$, $A \hi \dagger (h)$ is the unique member $g$ of $\cran (A)$  {satisfying} $A(g) = h$. We call $A \hi \dagger$ the Moore-Penrose inverse of $A$. Note that, according to our definition, $A \hi \dagger$ may or may not be bounded. An unbounded linear operator is essentially unestimable, because it is discontinuous. Nevertheless, $A \hi \dagger$ will never appear alone in our discussion: we see $A \hi \dagger$ only in the form $A \hi \dagger B$ where $B$ is another operator, say from $\ca K \to \ca H$. As discussed in \cite{li2017}, it is often reasonable to impose boundedness, compactness, or other similar assumptions on $A \hi \dagger B$, even when $A \hi \dagger$ itself is unbounded.  In the context of our applications, $A$ is usually a compact or trace-class operator, in which case, unless $A$ has finite rank, $A \hi \dagger$ is an unbounded operator. This  means $A \hi \dagger$ is unbounded unless we are in a finite-dimensional setting.

The following property of the Moore-Penrose inverse is useful. The proof is essentially that of Theorem 3.5.8 of \cite{Hsing2015}, though our definition of the Moore-Penrose inverse is slightly different from theirs.

\begin{proposition}\label{proposition:mp inverse} If $A \hi \dagger: \ran(A) \to \cran(A)$ is the Moore-Penrose inverse of $A$, then $A \hi \dagger A$ is the projection operator on to $\cran(A)$.
\end{proposition}

\subsection{Nonlinear sufficient dimension reduction}

Let $(\Omega \lo Y, \ca F \lo Y)$ be a measurable space, and let  $\ca H \lo Y$   be a Hilbert space  of functions defined on   $\Omega \lo Y$ with   $\ca H \lo Y \subseteq L \lo 2 (P \lo Y)$.
In the following, for a class $\ca S$ of functions  defined on $\Omega \lo X$ that are measurable with respect to $\ca F \lo X$, we use
\begin{align*}
\sigma \{ f (X): f \in  {\ca S} \}
\end{align*}
to denote the smallest $\sigma$-field that makes every $f(X)$, $f \in   {\ca S}$, measurable.
We now give a general formulation of nonlinear sufficient dimension reduction.

\begin{assumption}\label{assumption:nonlinear sdr} There exists a subset  $\ca A$ of $\ca H \lo X$  such that
\begin{align}\label{eq:sdr assumption}
Y \indep X | \sigma   \{ f(X): f \in \ca A \}.
\end{align}
Furthermore, $ \sigma   \{ f(X): f \in \ca A \}$ is the smallest $\sigma$-field such that $X$ and $Y$ are conditionally independent given it.
\end{assumption}
By requiring $ \sigma   \{ f(X): f \in \ca A \}  $ to be the smallest $\sigma$-field that makes (\ref{eq:sdr assumption}) holds, we are in effect assuming there is no redundant function in the set $\{ f(X): f \in \ca A \} $. This $\sigma$-field is called the central $\sigma$-field, denoted by $\ca G \lo {Y|X}$, and $ {(\ca H \lo X )} \lo {\ca G \lo {Y|X}}$ is called the central class, denoted by $\frak S \lo {Y|X}$. Henceforth,  we will abbreviate expressions such as (\ref{eq:sdr assumption}) as
\begin{align*}
Y \indep X |  \{ f(X): f \in \ca A \}.
\end{align*}

Since constants are not important for conditional independence, We can, without loss of generality, assume that the central class is contained in the closure of the range of $\Sigma \lo {XX}$. This is formally shown in the next proposition.

\begin{proposition} The following statements hold:
\begin{enumerate}
\item If $\ca H \lo X$ does not contain any nonzero constant function, then $\cenclass \subseteq \cran (\Sigma \lo {XX})$; \vspace{-0.08in}
\item If $\ca H \lo X$ contains a nonzero function, and $\ca H \lo X \hii 0 = \cran (\Sigma \lo {XX})$, then
\begin{align*}
Y \indep X | \{ f(X): \, f \in ( \ca H \lo X \hii 0 ) \lo {\ca G \lo {Y|X}} \}.
\end{align*}
\end{enumerate}
\end{proposition}

\proof 1. Note that $\ker (\Sigma \lo {XX})$ consists of all constant functions. If $\ca H \lo X$ does not contain nonzero constant functions, then $\ker (\Sigma \lo {XX}) = \{0\}$. Then $\cran (\Sigma \lo {XX}) = \ker ( \Sigma \lo {XX} ) \hi {\perp \lo 1} = \ca H \lo X$. Hence $\cenclass \subseteq \cran (\Sigma \lo {XX})$.

2. Since $\ker (\Sigma \lo {YX})$ is the class of constant functions, its orthogonal complement  $\cran(\Sigma \lo {XX})$ is the set
\begin{align*}
\left\{ f - \frac{\langle f , 1 \rangle \lo {\ca H \lo X}} { \langle 1, 1 \rangle \lo {\ca H \lo X}}: f \in \ca H \lo X \right\}.
\end{align*}
So, for each $f \in \ca A$,
\begin{align*}
\tilde f   = f  -\frac{\langle f   , 1 \rangle \lo {\ca H \lo X}} { \langle 1, 1 \rangle \lo {\ca H \lo X}}
\end{align*}
is a member  of $\cran(\Sigma \lo {XX})$. Since $\{f(X): f \in \ca A\}$  and $\{ \tilde f (X) : f \in \ca A \}$ generate the same $\sigma$-field, we have
\begin{align*}
Y \indep X | \{ \tilde f   (X): f \in \ca A \}.
\end{align*}
The asserted statement holds because $\{ \tilde f   (X): f \in \ca A \}$ and  $\{ f(X): \, f \in ( \ca H \lo X \hii 0 ) \lo {\ca G \lo {Y|X}} \}$ generate the same $\sigma$-field. \eop

This proposition shows that we can, without loss of generality assume that $\cenclass \subseteq \cran(\Sigma \lo {XX})$. We make this formal assumption below.

\begin{assumption}\label{assumption:in range}
\quad $\cenclass \subseteq \cran(\Sigma \lo {XX})$.
\end{assumption}

\def\cenclass{\frak{S} \lo {Y|X}}
\def\lat{\mathrm{Lat}}

The goal of nonlinear sufficient dimension reduction is to estimate the central class $\frak{S} \lo {Y|X}$. This usually proceeds as follows. Let $\frak{F}$ be the class of all distributions of $(X, Y)$. Let $\lat (\ca H \lo X)$ be the class of all closed linear subspaces of  $\ca H \lo X$. Here, the symbol Lat represents the word ``lattice'', because
$\lat (\frak{S} \lo {Y|X})$ is indeed a lattice in terms of the operations
\begin{align*}
\ca S \lo 1 \wedge \ca S \lo 2 = \ca S \lo 1 \cap \ca S \lo 2, \quad \ca S \lo 1 \vee \ca S \lo 2 = \overline{\spn} ( \ca S \lo 1 + \ca S \lo 2),
\end{align*}
where $ \ca S \lo 1 + \ca S \lo 2$ the set $\{a + b: a \in \ca S \lo 1, b \in \ca S \lo 2\}$.  Let $F \lo 0$ be the true distribution of $(X,Y)$ and $F \lo n$ the empirical distribution of $(X,Y)$ based on an i.i.d. sample $(X \lo 1, Y \lo 1), \ldots, (X \lo n , Y \lo n)$. Let $T: \frak{F} \to \lat (\ca H \lo X )$ be a mapping that sends a distribution in $\frak F$ to a closed subspace of $\ca H \lo X$. The lattice $\lat ( \ca H \lo X)$ is the parameter space for nonlinear sufficient dimension reduction, and the central class $\cenclass$ is the true parameter to be estimated. The mapping $T$ is called a statistical functional; $T (F \lo n)$ is the estimator, and $T ( F \lo 0)$ is usually the parameter value to which $T (F \lo n)$ converges. We now give a formal definition of the unbiasedness, exhaustiveness, and Fisher consistency of $T(F \lo n)$ as an estimator of the central class $\cenclass$.

\begin{definition}
    We say that an estimate $T (F \lo n)$ is unbiased for $\cenclass$ if
    $
        \cran ( R \lo {XY} ) \subseteq \cenclass,
   $
    exhaustive if
   $
        \cran ( R \lo {XY} ) \supseteq \cenclass,
     $
    and Fisher consistent if both hold.
\end{definition}

\subsection{Regression operator and GSIR}

We make the following assumption about $\ca H \lo Y$ and $L \lo 2 (P \lo Y)$, which is parallel  to Assumption \ref{assumption:norm norm} about $\ca H \lo X$ and $L \lo 2 (P \lo X)$.

\begin{assumption}\label{assumption:norm less norm Y} There exists a constant  $C > 0$  such that, for any $f \in \ca H \lo Y$,
\begin{align*}
  \| f \| \lo {L \lo 2 (P \lo  {Y})} \le C \| f \| \lo {\ca H \lo Y}.
\end{align*}
\end{assumption}

Consider the linear functionals
\begin{align*}
\ali T \lo 1: \ca H \lo X \to \real, \quad T \lo 1 (f) = E f(X),  \\
\ali T \lo 2: \ca H \lo Y \to \real, \quad T \lo 2 (g) = E g(Y),
\end{align*}
and the bilinear forms
\begin{align*}
\ali b \lo 1: \ca H \lo X \times \ca H \lo X \to \real, \quad b \lo 1 (f,g) = \cov[f(X), g(X)], \\
\ali b \lo 2: \ca H \lo X \times \ca H \lo Y \to \real, \quad b \lo 3 (f,g) = \cov[f(X), g(Y)], \\
\ali b \lo 3: \ca H \lo Y \times \ca H \lo X \to \real, \quad b \lo 3 (f,g) = \cov[f(Y), g(X)],  \\
\ali b \lo 4: \ca H \lo Y \times \ca H \lo Y \to \real, \quad b \lo 2 (f,g) = \cov[f(Y), g(Y)].
\end{align*}
It can be easily shown that, under Assumptions \ref{assumption:norm norm} and \ref{assumption:norm less norm Y}, these functionals and bilinear forms are bounded. We record these facts below without proof.
\begin{lemma} Under Assumptions \ref{assumption:norm norm} and \ref{assumption:norm less norm Y}, the linear functionals $T \lo 1$ and $T \lo 2$ are bounded, and the bilinear forms $b \lo 1, b \lo 2, b \lo 3,  b \lo 4$ are bounded.
\end{lemma}

By Riesz representation theorem, there exist $\mu \lo X \in \ca H \lo X$ and $\mu \lo Y \in \ca H \lo Y$ such that
\begin{align*}
T \lo 1 (f) = \ali \langle f, \mu \lo X \rangle \lo {\ca H \lo X} \ \mbox{for all $f \in \ca H \lo X$};  \\
T \lo 2 (g) = \ali \langle g, \mu \lo Y \rangle \lo {\ca H \lo Y} \ \mbox{for all $f \in \ca H \lo Y$}.
\end{align*}
We call $\mu \lo X$ and $\mu \lo Y$ the mean elements in  $\ca H \lo X$ and  $\ca H \lo Y$, respectively. Furthermore, by Theorem 2.2 of \cite{Conway1990}, there exist operators
\begin{align*}
 \Sigma \lo {XX} \in \ca B ( \ca H \lo X, \ca H \lo X) , \
 \Sigma \lo {XY} \in \ca B ( \ca H \lo X, \ca H \lo Y) , \
 \Sigma \lo {YX} \in \ca B ( \ca H \lo Y, \ca H \lo X) , \
 \Sigma \lo {YY} \in \ca B ( \ca H \lo Y, \ca H \lo Y)
\end{align*}
such that
\begin{align*}
b \lo 1 (f, g) = \ali \cov[f(X), g(X)] =   \langle f, \Sigma \lo {XX} g \rangle \lo {\ca H \lo X} \quad \mbox{for all  $f, g \in \ca H \lo X$},\\
b \lo 2 (f, g) = \ali \cov[f(X), g(Y)] =   \langle f, \Sigma \lo {XY} g \rangle \lo {\ca H \lo X} \quad \mbox{for all  $f \in \ca H \lo X, g \in \ca H \lo Y$},\\
b \lo 3 (f, g) = \ali \cov[f(Y), g(X)] =   \langle \Sigma \lo {YX}f,  g \rangle \lo {\ca H \lo Y} \quad \mbox{for all  $f \in \ca H \lo Y, g \in \ca H \lo X$},\\
b \lo 4 (f, g) = \ali \cov[f(Y), g(Y)] =   \langle f, \Sigma \lo {YY} g \rangle \lo {\ca H \lo Y} \quad \mbox{for all  $f, g \in \ca H \lo Y$}.
\end{align*}
The operator $\Sigma \lo {XX}$ is called the covariance operator in $\ca H \lo X$, $\Sigma \lo {XY}$ the covariance operator from $\ca H \lo Y$ to $\ca H \lo X$, $\Sigma \lo {YX}$ the covariance operator from $\ca H \lo X$ to $\ca H \lo Y$, and $\Sigma \lo {YY}$ the covariance operator from $\ca H \lo Y$ to $\ca H \lo Y$.

We next introduce the regression operator. To do so we make the following assumption.
\begin{assumption}\label{assumption:ran in ran} $\ran (\Sigma \lo {XY} ) \subseteq \ran ( \Sigma \lo {XX} )$.
\end{assumption}
As argued in \cite{li2017}, this assumption is about the smoothness in the relation between $X$ and $Y$.
Under Assumption \ref{assumption:ran in ran}, the linear operator
\begin{align*}
R \lo {XY}: \ca H \lo Y \to \ca H \lo X, \quad R \lo {XY} = \Sigma \lo {XX} \hi \dagger \Sigma \lo {XY}
\end{align*}
is well defined. We call this operator the regression operator. We make the following assumptions about the operators $\Sigma \lo {XX}$, $\Sigma \lo {YY}$, and $R \lo {XY}$.

\begin{assumption}\label{assumption:compact} The operators $\Sigma \lo {XX}$ and  $\Sigma \lo {YY}$ are compact.
\end{assumption}

This assumption is   very mild. In fact, if $\ca H \lo X$ and $\ca H \lo Y$ are RKHS's, then it is well known that $\Sigma \lo {XX}$ and $\Sigma \lo {YY}$ trace class operators, and therefore compact.

\begin{assumption}\label{assumpption:compact for R}
\quad  $R \lo {XY}$ are compact operators.
\end{assumption}

Again, as  argued in \cite{li2018sufficient}, requiring $R \lo {XY}$ to be compact amounts to  imposing a degree of smoothness on the relation between $X$ and $Y$.

Consider any statistical functional that satisfies the condition
\begin{align}\label{eq:condition for gsir}
T(F \lo 0) = \cran ( R \lo {XY} ),
\end{align}
where the right-hand side is the regression operator based on the true distribution of $(X,Y)$. We now give a formal definition of the genearlized sliced inverse regression, or GSIR. See \cite{lee2013} and \cite{li2018sufficient}.

\begin{definition}
 Any statistical functional that satisfies (\ref{eq:condition for gsir}) called the generalized sliced inverse regression, or GSIR.
\end{definition}
The motivation for calling this estimator the generalized sliced inverse regression is that it resembles sliced inverse regression (SIR) of \cite{Li1991}: if we replace the scalar product $\beta \trans X$ in the eigenvalue problem that defines SIR  by the RKHS inner product $\langle f, \ka \lo X (\cdot, X) \rangle \lo {\ca H \lo X}$, then we obtain GSIR. See \cite{li2018sufficient}, page 215.


\section{Unbiasedness and Fisher consistency of GSIR}\label{section:unbiased}

In this section, we prove the unbiasedness and Fisher consistency of the closure of the range of the regression operator  using the new definition of relative universality. Towards the end of this section we will also discuss the gap in \cite{li2018sufficient}'s proof. We begin with unbiasedness.

\subsection{Unbiasedness}

We first prove a lemma, which gives an equivalent condition for a function to be a member of  $[L \lo 2 (P \lo X) \lo { \ca G}] \hi {\perp \lo 3}$.

\begin{lemma}\label{lemma:constant} Suppose $\ca G$ is a sub-$\sigma$-field of $\ca F $. Then $f \in [L \lo 2 (P \lo X) \lo { \ca G}] \hi {\perp \lo 3}$ if and only if
\begin{align}\label{eq:E and E}
E[f(X) | \ca G] = E[f(X)] \quad \mbox{almost surely}.
\end{align}
\end{lemma}

\proof Let  $f \in [L \lo 2 (P \lo X) \lo { \ca G}] \hi {\perp \lo 3}$ and $g \in L \lo 2 (P \lo X) \lo { \ca G}$. Then
\begin{align*}
f \in L \lo 2 (P \lo X) \ \mbox{and} \ \cov[E(f(X)|\ca G), g(X)] = \cov[f(X), g(X)] = 0.
\end{align*}
In particular, taking  $g(X)= E ( f(X)| \ca G )$, we have
$
\var [  E ( f(X)| \ca G)] = 0,
$
which implies $E(f(X) |\ca G) = \mbox{constant}$ almost surely. Taking unconditional expectation on both sides, we have the second relation in (\ref{eq:E and E}). The first relation holds because $L \lo 2 (P \lo X) \lo { \ca G} \subseteq L \lo 2 (P \lo X)$.

Suppose $f$ satisfies (\ref{eq:E and E}) and $g   \in L \lo 2 (P \lo X) \lo { \ca G}$. Then
\begin{align*}
\cov[f(X), g(X)] =   \cov[E(f(X)|\ca G), g(X)] =   \cov[E(f(X)), g(X)] = 0.
\end{align*}
Hence $f \in  [L \lo 2 (P \lo X) \lo { \ca G}] \hi {\perp \lo 3}$. \eop

We are now ready to prove the unbiasedness of GSIR.

\begin{theorem}\label{theorem:unbiasedness} Suppose Assumptions \ref{assumption:norm norm} through \ref{assumpption:compact for R} are satisfied. If $\ca H \lo X$ {is dense in $L \lo 2 (P \lo X)$ modulo constants}, then
\begin{align}\label{eq:unbiased}
\cran ( R \lo {XY} ) \subseteq \frak S \lo {Y|X}.
\end{align}
\end{theorem}

\proof We first show that
\begin{align}\label{eq:first show}
\cran( \Sigma \lo {XY})  \subseteq \Sigma \lo {XX} \frak S \lo {Y|X}.
\end{align}
 Since
\begin{align*}
\frak S \lo {Y|X} =( \ca H \lo X ) \lo {{\ca G} \lo {Y|X}} = \cap \lo {\epsilon > 0} ( \ca H \lo X ) \lo {{\ca G} \lo {Y|X}} (\epsilon),
\end{align*}
it suffices to show that
$
\ran( \Sigma \lo {XY})  \subseteq  \Sigma \lo {XX} (\ca H \lo X ) \lo {{\ca G} \lo {Y|X}} (\epsilon)
$
for any $\epsilon > 0$. Or equivalently, for any $\epsilon > 0$,
\begin{align*}
[  \Sigma \lo {XX} (\ca H \lo X ) \lo {{\ca G} \lo {Y|X}} (\epsilon) ] \hi {\perp \lo 1}\subseteq [ \cran( \Sigma \lo {XY}) ]  \hi {\perp \lo 1} = \ker ( \Sigma \lo {YX}).
\end{align*}
Let $f \in [  \Sigma \lo {XX} (\ca H \lo X ) \lo {{\ca G} \lo {Y|X}} (\epsilon) ] \hi {\perp \lo 1}$. Since $[  \Sigma \lo {XX} (\ca H \lo X ) \lo {{\ca G} \lo {Y|X}} (\epsilon) ] \hi {\perp \lo 1}\subseteq [(\ca H \lo X ) \lo {{\ca G} \lo {Y|X}} (\epsilon) ] \hi {\perp \lo 3}$, we have  $f \in[(\ca H \lo X ) \lo {{\ca G} \lo {Y|X}} (\epsilon) ] \hi {\perp \lo 3}$.  {Since $\ca H \lo X$ is dense in $L \lo 2 (P \lo X)$ modulo constants, by Theorem \ref{theorem:relative universality}, $\ca H \lo X$ is relative universal with respect to $\ca G \lo {Y|X}$}.
Hence $f \in [ L \lo 2 (P \lo X) \lo {\ca G \lo {Y|X}} ]  \hi {\perp \lo 3}$. By Lemma \ref{lemma:constant},
\begin{align*}
E[f(X) | \ca G \lo {Y|X} ] = E [f (X)].
\end{align*}
Since $\ca G \lo {Y|X}$ is sufficient,
\begin{align*}
E[f(X) | Y] = E[E(f(X)|Y, \ca G \lo {Y|X}) | Y]  = E[E(f(X)|\ca G \lo {Y|X}) | Y]  =  E [f (X)].
\end{align*}
So, for any $y \in \Omega \lo Y$,
\begin{align*}
(\Sigma \lo {YX} f)(y) = \ali  E [ ( \ka \lo Y (y, Y) - \mu \lo Y (y) ) \langle \ka \lo X (\cdot, X)  - \mu \lo X, f \rangle \lo {\ca H \lo X} ] \\
= \ali  E [ ( \ka \lo Y (y, Y) - \mu \lo Y (y) ) (f(X)-Ef(X)) ] \\
= \ali  E [ ( \ka \lo Y (y, Y) - \mu \lo Y (y) ) E(f(X)-Ef(X)|Y) ] =0,
\end{align*}
which proves (\ref{eq:first show}).

Next, applying $\Sigma \lo {XX} \hi \dagger$ on the left of the both sides of the equation (\ref{eq:first show}), we have
\begin{align}\label{eq:mp inverse cran}
\Sigma \lo {XX} \hi \dagger\cran( \Sigma \lo {XY})  \subseteq \Sigma \lo {XX} \hi \dagger\Sigma \lo {XX} \frak S \lo {Y|X}.
\end{align}
By Proposition \ref{proposition:mp inverse},  $\Sigma \lo {XX} \hi \dagger \Sigma \lo {XX} $ is the projection on to $\cran ( \Sigma \lo {XX} )$, which, together with Assumption \ref{assumption:in range}, implies that  $\Sigma \lo {XX} \hi \dagger\Sigma \lo {XX} \frak S \lo {Y|X}= \frak S \lo {Y|X}$. The left-hand side of (\ref{eq:mp inverse cran}) can be rewritten as $\cran( \Sigma \lo {XX} \hi \dagger\Sigma \lo {XY}) = \cran (R \lo {XY})$. Hence (\ref{eq:unbiased}) holds.
\eop


\vspace{-.2in}

\subsection{Exhaustiveness and Fisher consistency}

We now turn to exhaustiveness and Fisher consistency. As in \cite{lee2013}, we say that a sub $\sigma$-field $\ca G$ of $\ca F$ is complete if, for any $f \in L \lo 2 (P \lo X) \lo {\ca G}$,
\begin{align*}
E[f (X) | Y] = \mbox{constant} \ \mbox{almost surely} \ \Rightarrow \ f(X) = \mbox{constant} \ \mbox{almost surely}.
\end{align*}
The next theorem gives the sufficient condition for exhaustiveness. This result has not been recorded in the literature previously, though the proof follows easily from that of Theorem 13.2 of Li (2018), and is therefore omitted.

\begin{theorem}\label{theorem:exhaustiveness} Suppose Assumptions \ref{assumption:norm norm} trough \ref{assumpption:compact for R} are satisfied. Furthermore, suppose
\begin{enumerate}
\item $\ca H \lo Y$ is dense in $L \lo 2 (P \lo Y)$ modulo constants;
\vspace{-0.08in}
\item $\ca G \lo {Y|X}$ is complete.
\end{enumerate}
Then
$
\cran ( R \lo {XY} ) \supseteq \frak S \lo {Y|X}.
 $
\end{theorem}

Interestingly, for exhaustiveness, we do not require $\ca H \lo X$ to be dense in $L \lo 2 (P \lo X)$ modulo constants.
Combining Theorem \ref{theorem:unbiasedness} and Theorem \ref{theorem:exhaustiveness}, we arrive at the following result, which is essentially Theorem 13.2 and Theorem 13.3 of \cite{li2018sufficient} combined,  though, as mentioned before, we do not require $\ca H \lo X$ or $\ca H \lo Y$ to be RKHS.

\begin{theorem} Suppose Assumptions \ref{assumption:norm norm} through \ref{assumpption:compact for R} are satisfied, and
\begin{enumerate}
\item $\ca H \lo X$ is dense in $L \lo 2 (P \lo X)$ modulo constants;
\vspace{-0.08in}
\item $\ca H \lo Y$ is dense in $L \lo 2 (P \lo Y)$ modulo constants;
\vspace{-0.08in}
\item $\ca G \lo {Y|X}$ is complete.
\end{enumerate}
Then
$
\cran ( R \lo {XY} ) = \frak S \lo {Y|X}.
 $
\end{theorem}

\subsection{The gap in \cite{li2018sufficient}'s proof}

To provide more backgrounds and insights into the development of this paper,   we now give a detailed description of the gap  in the proof of Theorem 13.3 of \cite{li2018sufficient},
In our more general setting, Theorem 13.3 in \cite{li2018sufficient} can be stated as follows:
\begin{align}\label{quote:li2018}
\begin{split}
 \ali \mbox{\em For any  given sub-$\sigma$-field   $\ca G$  of $\sigma (X)$. If $\ca H \lo X$ is dense  in $L \lo 2 (P \lo X)$ modulo} \\
\ali   \mbox{\em constants,  then $(\ca H \lo X) \lo {\ca G} $ is dense in $L \lo 2 (P \lo X) \lo {\ca G}$ modulo constants. }
    \end{split}
\end{align}
Intuitively, the statement says that if $\ca H \lo X$ is rich enough to approximate any member of $L \lo 2 (P \lo X)$, then $(\ca H \lo X) \lo {\ca G}$ is rich enough to approximate any member of $L \lo 2 (P \lo X) \lo {\ca G}$, which seems to be  a plausible statement.

Let $\sim$ be the equivalent relation in the proof of Lemma \ref{lemma:perp 3}, and let $L \lo 2 (P \lo X ) / \sim$ be the quotient space with respect to $\sim$. It can be easily shown that this  quotient space is a Hilbert space in terms of the inner product $\langle f, g \rangle = \cov [ f(X), g(X)]$. In the following, when we say a function $f$ is a member of $L \lo 2 (P \lo X)/ \sim$,  we mean the equivalence class $\{f + c: c \in \real \}$ is a member of $L \lo 2(P \lo X) / \sim$.  The proof in \cite{li2018sufficient}   proceeds roughly in following five steps (the details can be found on page 216 of \cite{li2018sufficient}).
\begin{enumerate}
    \item Let $\ca A = \{h \lo n: n = 1, 2, \cdots \}$ be the set of eigenfunctions of  $\Sigma \lo {XX}$. Group them into subsets $\ca A \lo {\ca G}$, consisting   those eigenfunctions   that are measurable   $\ca G $, and $\ca A \lo {\ca G} \hi c$, consisting of those that are not.
    \item Let $\ca M \lo {\ca G}$ be the closure of $\ca A \lo {\ca G}$ in $L \lo 2 (P \lo X)/\sim$, and $\ca M \lo {\ca G} \hi {\perp \lo 3}$ the closure of $\ca A \setminus \ca A \lo {\ca G}$ in $L \lo 2 (P \lo X ) / \sim$.
    \item Let $f$ be a member of $L \lo 2 (P \lo X) \lo {\ca G}$ and let $\{s \lo n \}$ be a sequence in $\ca H \lo X$ such that $\var[s \lo n (X) - f(X)] \to 0$;  this is possible because $\ca H \lo X$ is dense in $L \lo 2 (P \lo X)$ modulo constants.
    \item Decompose $s \lo n$ as $s \lo n \hii 1 + s \lo n \hii 2$, where  $s \lo n \hii 1 \in \ca M \lo {\ca G}$ and $s \lo n \hii 2 \in {\ca M } \lo {\ca G} \hi {\perp \lo 3}$ (this is possible because $\ca H \lo X \subseteq \ca M \lo {\ca G} + \ca M \lo {\ca G} \hi {\perp \lo 3}$), and show that they are Cauchy sequences in $L \lo 2 (P \lo X ) / \sim$. Let $s \hii 1$ and $s \hii 2$ be the limit of $s \lo n \hii 1$ in $\ca M \lo {\ca G}$ and $s \hii 2$ the limit of $s \lo n \hii 2$ in ${\ca M} \lo {\ca G} \hi {\perp \lo 3}$.
    \item Show that $s \hii 2 = 0$ and whence that $s \lo n \hii 1$ converges to $f$ in $L \lo 2 (P \lo X) / \sim$. Conclude that $(\ca H \lo X)\lo {\ca G}$ is dense in $L \lo 2 (P \lo X) \lo {\ca G}$ modulo constants.
\end{enumerate}
The problem with this proof is that $\ca M \lo {\ca G}$ is the $L \lo 2 (P \lo X) / \sim$ closure of $\ca A$, rather than  the $\ca H \lo X$ closure of $\ca A$. Thus the sequence $s \lo n \hii 1$ need not be members of $(\ca H \lo X ) \lo {\ca G}$. Thus we have only shown that there is a sequence in the $L \lo 2 (P \lo X) / \sim$ closure of $(\ca H \lo X) \lo {\ca G}$ that converges to $f$; we have not shown that there is a sequence in $(\ca H \lo X) \lo {\ca G}$ that converges to $f$. This is  the gap.
What we did in the new proof is to replace  measurable with respect to $\ca G$ by $\epsilon$-measurable with respect to $\ca G$ to get around the  problem.

\section{Conclusion}\label{section:conclusion}

In this paper we rigorously define  the notion of relative universality and crystalize   its role in characterizing conditional independence. That is, through relative universality, we established that the range of the regression operator generates the sub $\sigma$-field of $\ca F$ given  which $Y$ and $X$ are conditionally  independent. More specifically, our key result is this:
\begin{quote}
   {\em   If the regression operator $R \lo {XY}$ is defined and compact, $\ca H \lo X$ is dense in $L \lo 2 (P \lo X)$ modulo constants, and $\ca G \lo {Y|X}$ is the smallest $\sigma$-field such that $Y \indep X | \ca G \lo {Y|X}$, then
    \begin{align*}
        \sigma \{f(X): f \in \cran(R \lo {XY} )\} \subseteq \ca G \lo {Y|X}.
    \end{align*}
    If, furthermore, $\ca G \lo {Y|X}$ is complete, then the equality holds. }
\end{quote}
This result precisely describes the relation between the regression operator and conditional independence. The significance of this result is that the regression operator can be estimated by replacing the moments in it with sample averages \citep{li2018sufficient}.  We proved this via relative universality, a modified version of this concept in \cite{li2018sufficient}.

To summarize the logic line that leads to the modified definition of relative universality,   consider the following three statements:
\begin{enumerate}
    \item $\ca H \lo X$ is dense in $L \lo 2 (P \lo X)$ modulo constants;
    \vspace{-.08in}
    \item the set of functions in $\ca H \lo X$ that are measurable $\ca G$ is dense (modulo constants) in the set of functions in $L \lo 2 (P \lo X)$ that are measurable $\ca G$;  \vspace{-.08in}
    \item for each $\epsilon> 0$,
    the set of functions in $\ca H \lo X$ that are $\epsilon$-measurable with respect to $\ca G$ is dense (modulo constants) in the set of functions in $L \lo 2 (P \lo X)$ that are measurable with respect to $\ca G$.
\end{enumerate}
\cite{li2018sufficient} attempted to prove the unbiasedness and Fisher consistency of $\cran (R \lo {XY})$ using  the assertion that 1 implies 2 (Theorem 13.3), but this assertion may not be right --- at least there is a gap in the proof Theorem 13.3. Our new proof uses the fact that 1 implies 3, which is rigorously proved in this paper.

While the notion of relative universality was originally introduced in the context of nonlinear sufficient dimension reduction, the significance of the current note is beyond this context.  Indeed, it is a general  mechanism through which we establish that the regression operator characterizes conditional independence. Since conditional independence is widely used in statistics,  machine learning, and many other scientific disciplines,   the theoretical framework rigorously established in this paper will be  useful for developing and    applying this important methodology.

In concluding this paper, we should mention that there is still a possibility that statement (\ref{quote:li2018}) may turn out to be  correct, for, as we mentioned, it sounds reasonable. We will leave it as an open problem, to be resolved either by finding a reference that we are unaware of or by proving it via another route from that used in  \cite{li2018sufficient}. Regardless of the correctness of (\ref{quote:li2018}), though, we now have  rigorously established the relation between conditional independence and regression operator via the modified form of relative universality.

\vskip 0.2in
\bibliographystyle{unsrtnat}
\bibliography{references, proposal}

\end{document}